\documentclass[12pt,english]{article}
\usepackage{theorem}
\usepackage{a4wide}
\usepackage{graphics}
\usepackage{epsfig}
\usepackage{amssymb}
\usepackage{latexsym}
\usepackage{amsmath}

\newtheorem{thh}{Theorem}[section]
\newtheorem{df}[thh]{Definition}

\newtheorem{lem}[thh]{Lemma}
\newtheorem{cor}[thh]{Corollary}
\newtheorem{prop}[thh]{Proposition}

\title{Minimal invariant varieties and first integrals \\ for algebraic foliations}
\author{Philippe Bonnet}
\date{}

\newcommand{\dem}{{\em Proof: }}
\newcommand{\qed}{\begin{flushright} $\blacksquare$\end{flushright}}
\newcommand{\der}{ \partial}

\newcommand{\FF}{{\cal{F}}}
\newcommand{\CC}{\mathbb C}

\newcommand{\OX}{{\cal{O}}_{X}}

\begin{document}
\maketitle

\begin{center} { \small
Mathematisches Institut, Universit\"at Basel\\
Rheinsprung 21, Basel 4051, Switzerland\\
e-mail: bonnet@math-lab.unibas.ch}
\end{center}

\begin{abstract}
Let $X$ be an irreducible algebraic variety over $\CC$, endowed with an algebraic foliation $\FF$.
In this paper, we introduce the notion of minimal invariant variety $V(\FF,Y)$ with respect to $(\FF,Y)$, where $Y$
is a subvariety of $X$. If $Y=\{x\}$ is a smooth point where the foliation is regular, its minimal invariant
variety is simply the Zariski closure of the leaf passing through $x$. First we prove that for very generic
$x$, the varieties $V(\FF,x)$ have the same dimension $p$. Second we generalize a result due to X.Gomez-Mont (see \cite{G-M}). More
precisely, we prove the existence of a dominant rational map $F:X\rightarrow Z$, where $Z$ has dimension $(n-p)$, such that
for very generic $x$, the Zariski closure of $F^{-1}(F(x))$ is one and only one minimal invariant variety of a point.
We end up with an example illustrating both results. 
\end{abstract}

\section{Introduction}

Let $X$ be an affine irreducible variety over $\CC$, and $\OX$ its ring of regular functions. Let
$\FF$ be an algebraic foliation, i.e. a collection of algebraic vector fields on $X$ stable by Lie
bracket. We consider the elements of $\FF$ as $\CC$-derivations on the ring $\OX$. In this paper,
we are going to extend the notion of algebraic solution for $\FF$: this will be the minimal tangent
varieties for $\FF$. We will study some of their properties and relate them to the existence of rational
first integrals for $\FF$. 

Recall that a subvariety $Y$ of $X$ is an algebraic solution of $\FF$ if $Y$ is the closure (for the metric topology)
of a leaf of $\FF$. A non-constant rational function $f$ on $X$ is a first integral if $\der(f)=0$ for any $\der$ in $\FF$. Since the works
of Darboux, the existence of such varieties has been extensively
studied in the case of codimension 1 foliations (see \cite{Jou},\cite{Gh},\cite{Bru}). In particular,
from these works, we know that only two cases may occur for codimension 1 foliations:

\begin{itemize}
\item{$\FF$ has finitely many algebraic solutions,}
\item{$\FF$ has infinitely many algebraic solutions, and a rational first integral.}
\end{itemize}
So rational first integrals appear if and only if all leaves of $\FF$ are algebraic solutions. In this case, the fibres of any
rational first integral is a finite union of closures of leaves. This fact has been generalised by Gomez-Mont
(see \cite{G-M}) in the following way.

\begin{thh}
Let $X$ be a projective variety and $\FF$ an algebraic foliation on $X$ such that all leaves are quasi-projective.
Then there exists a rational map $F:X\to Y$ such that, for every generic point $y$ of $Y$, the Zariski closure
of $F^{-1}(y)$ is the closure of a leaf of $\FF$.
\end{thh}
We would like to find a version of this result that does not need all leaves to be algebraic. To that purpose, we need to give a correct
definition to the algebraic object closest to a leaf. A good candidate would be the Zariski closure of a leaf,
but this choice may rise difficulties due to the singularities of both $X$ and $\FF$. We counterpass this problem by the following
algebraic approach.

Let $Y$ be an algebraic subvariety of $X$ and $I_Y$ the ideal of vanishing functions on $Y$. Let ${\cal{J}}$ be the set of ideals
$I$ in $\OX$ satisfying the two conditions:
$$
(i) \quad (0) \; \subseteq I \; \subseteq I_Y \quad \quad {\rm{and}} \quad \quad (ii) \quad \forall \; \der \; \in \; \FF, \quad \der(I)\; \subseteq\; I
$$
Since $(0)$ belongs to ${\cal{J}}$, ${\cal{J}}$ is non-empty and it is partially ordered by the inclusion. Since
it is obviously inductive, ${\cal{J}}$ admits a maximal element $I$. If $J$ is any other ideal of ${\cal{J}}$,
then $I+J$ enjoys the conditions $(i)$ and $(ii)$, hence it belongs to ${\cal{J}}$. By maximality, we have
$I=I+J$ and $J$ is contained in $I$. Therefore $I$ is the unique maximal element of ${\cal{J}}$, which
we denote by $I(\FF,Y)$. 

\begin{df}
The minimal invariant variety $V(\FF,Y)$ is the zero set of $I(\FF,Y)$ in $X$.
\end{df}
From a geometric viewpoint, $V(\FF,Y)$ can be seen as the smallest subvariety containing $Y$ and invariant by the flows of
all elements of $\FF$. In particular, if $x$ is a smooth point of $X$ where the foliation is regular, then $V(\FF,x)$ is
the Zariski closure of the leaf passing through $x$. In section \ref{order}, we show that $V(\FF,Y)$ is irreducible if $Y$
is itself irreducible. \\

In this paper, we would like to study the behaviour of these invariant varieties, and relate it to the
existence of first integrals. We analyze some properties of the function:
$$
n_{\FF} : X \longrightarrow \mathbb{N}, \quad x \longmapsto dim\; V(\FF,x)
$$
Let ${\cal{M}}$ be the $\sigma$-algebra generated by the Zariski topology on $X$. A function $f: X\rightarrow
\mathbb{N}$ is {\em measurable for the Zariski topology} if $f^{-1}(p)$ belongs to ${\cal{M}}$ for any $p$.
The space ${\cal{M}}$
contains in particular all countable intersections $\theta$ of Zariski open sets. A property ${\cal{P}}$ holds for {\em
every very generic point $x$ in $X$} if ${\cal{P}}(x)$ is true for any point $x$ in such an intersection $\theta$.

\begin{thh} \label{mesure}
Let $X$ be an affine irreductible variety over $\CC$ and $\FF$ an algebraic foliation on $X$. Then the function $n_{\FF}$ is
measurable for the Zariski topology. Moreover there exists
an integer $p$ such that $(1)$
$n_{\FF}(x)\leq p$ for any point $x$ in $X$ and $(2)$ $n_{\FF}(x)=p$ for any very generic point $x$ in $X$. 
\end{thh}
Set $p=max\; dim \; V(\FF,x)$ and note that $p$ is achieved for every generic point of $X$. In the last section, we will
produce an example of a foliation on $\CC^4$ where the function $n_{\FF}$ is measurable but not constructible for the Zariski
topology. In this sense, theorem \ref{mesure} is the best result one can expect for any algebraic foliation.

Let $K_{\FF}$ be the
field generated by $\CC$ and the rational first integrals of $\FF$. By construction,
the invariant varieties $V(\FF,x)$ are defined set-theoretically, and they seem to appear randomly, i.e.
with no link within each other. In fact there does exist some order among them, and we are going to
see that they are "mostly" given as the fibres of a rational map. More precisely:

\begin{thh} \label{fibration}
Let $X$ be an affine irreducible variety over $\CC$ of dimension $n$ and $\FF$ an algebraic foliation on $X$.
Then there exists a dominant rational map $F: X\rightarrow Y$, where $Y$ is irreducible of dimension $(n-p)$,
such that for every very generic point $x$ of $X$, the Zariski closure of $F^{-1}(F(x))$ is equal to
$V(\FF,x)$. In particular,
the transcendence degree of $K_{\FF}$ over $\CC$ is equal to $(n-p)$. 
\end{thh}
The idea of the proof is to construct enough rational first integrals. These will
be the coordinate functions of the rational map $F$ given above. The construction
consists in choosing a codimension $d$ irreducible variety $H$ in $X$. We show
there exists an integer $r>0$ such that, for every very generic point $x$ of $X$,
$V(\FF,x)$ intersects $H$ in $r$ distinct
points $y_1,...,y_r$.
We then obtain a correspondence:
$$
{\cal{H}}: x\longmapsto  \{y_1,...,y_r\}
$$
We can modify ${\cal{H}}$ so as to get a rational map $F$ that represents every $r$-uple $\{y_1,...,y_r\}$
by a single point. Since the image of $x$ only depends on the intersection of $V(\FF,x)$ with $H$, the
map $F$ will be invariant with respect to the elements of $\FF$. \\

One question may arise after these two results. Does there exist an effective way of computing these
minimal invariant varieties and detect the presence of rational first integrals? For instance,
we may attempt to use the description of the ideals $I(\FF,Y)$ given by lemma \ref{autre}. Unfortunately we cannot hope
to compute them in a finite number of steps bounded, for instance, by the degrees of the components
of the vector fields of $\FF$. Indeed, consider the well-known derivation $\der$ on $\CC^2$:
$$
\der=px\frac{\der}{\der x} + qy\frac{\der}{\der y}
$$
For any couple of non-zero coprime integers $(p,q)$, this derivation will have $f(x,y)=x^qy^{-p}$ as a rational first integral,
and we cannot find another one of smaller degree. The minimal invariant varieties of points will be given in general by the fibres
of $f$. Therefore we cannot bound the degree of the generators of $I(\FF,x)$ solely by the degree of $\der$.

However, we may find them by an inductive process. For one derivation, an approach is given in the paper of J.V.Pereira
via the notion of extatic curves (see \cite{Pe}). The idea is to compute a series of Wronskians attached to the
derivation. Then one of them vanishes identically if and only the derivation has a rational first integral. 
\\

Last thing to say is that the previous results carry over all algebraic irreducible varieties. Given an algebraic
variety $X$ with an algebraic foliation, we choose a covering of $X$ by open affine sets $U_i$ and work
on the $U_i$. For any algebraic subvariety $Y$ of $X$, we define the minimal invariant variety $V(\FF,Y)$ by gluing
together the Zariski closure of the varieties $V(\FF,Y\cap U_i)$ in $X$.

\section{The contact order with respect to $\FF$} \label{order}

In this section, we are going to show that the minimal invariant variety $V(\FF,Y)$ is irreducible if $Y$ is
irreducible. This result is already known when $\FF$ consists of one derivation (see \cite{Ka}). We
could reproduce the proof given in \cite{Ka} for any set of derivations, but we prefer to adopt another
strategy. We will instead introduce a notion of contact order with respect to $\FF$, and we will use it to
show that $I(\FF,Y)$ is prime if $I_Y$ is prime. Denote by $M_{\FF}$ the $\OX$-module spanned by the
elements of $\FF$. We start by giving the following characterisation of $I(\FF,Y)$.

\begin{lem} \label{autre}
$\displaystyle I(\FF,Y)=\left\{ f \in I_Y, \; \forall \der_1,...,\der_k \in M_{\FF}, \; \der_1\circ ...\circ
\der_k(f)\in I_Y \right\}$
\end{lem}
\dem Let $f$ be an element of $I_Y$ such that $\der_1\circ ...\circ \der_k(f)$ belongs to $I_Y$ for any
$\der_1,...,\der_k$ in $M_{\FF}$. Then $\der_1\circ ...\circ \der_k(f)$ belongs to $I_Y$ for any elements
$\der_1,...,\der_k$ of $\FF$. Let $I$ be the ideal generated by $f$ and all the elements of the form
$\der_1\circ ...\circ \der_k(f)$, where every $\der_i$ lies in $\FF$. By construction, this ideal is
contained in $I_Y$, and is stable
by every derivation of $\FF$. Therefore $I$ is contained in $I(\FF,Y)$, and a fortiori $f$ belongs
to $I(\FF,Y)$. We then have the inclusion:
$$
\left\{ f \in I_Y, \; \forall \der_1,...,\der_k \in M_{\FF}, \; \der_1\circ ...\circ \der_k(f)\in I_Y \right\}\subseteq
I(\FF,Y)
$$
Conversely let $f$ be an element of $I(\FF,Y)$. Since $I(\FF,Y)$ is contained in $I_Y$ and is stable by
every derivation of $\FF$, $\der_1\circ ...\circ \der_k(f)$ belongs to $I_Y$ for any elements
$\der_1,...,\der_k$ of $\FF$. Since $M_{\FF}$ is spanned by $\FF$, $\der_1\circ ...\circ \der_k(f)$
belongs to $I_Y$ for any $\der_1,...,\der_k$ in $M_{\FF}$
\qed
Since the space
of $\CC$-derivations on $\OX$ is an $\OX$-module of finite type and $\OX$ is noetherian,
$M_{\FF}$ is finitely generated as an $\OX$-module. 
Let $\{\der_1,...,\der_r\}$ be a system of generators of $M_{\FF}$. If $I=(i_1,...,i_n)$ belongs to $\{1,...,r\}^n$,
we set $\der_I = \der_{i_1}\circ ...\circ \der_{i_n}$ and $|I|=n$. By convention $\{1,...,r\}^0=\{\emptyset\}$,
$|\emptyset|=0$ and
$\der_{\emptyset}$ is the identity on $\OX$. We introduce the following map:
$$
ord_{\FF,Y}: \OX \longrightarrow \mathbb{N} \cup \{+\infty\}, \quad f \longmapsto \inf\left\{|I|, \; \der_I(f)\not\in
I_Y\right\}
$$

\begin{df}
The map $ord_{\FF,Y}$ is the contact order with respect to $(\FF,Y)$.
\end{df}
By lemma \ref{autre}, $f$ belongs to $I(\FF,Y)$ if and only if $ord_{\FF,Y}(f)=+\infty$, and $f$ does not belong to
$I_Y$ if and only if $ord_{\FF,Y}(f)=0$. A priori, the map $ord_{\FF,Y}$
depends on the set of generators chosen for $M_{\FF}$. We are going to see that it only depends on $\FF$. Let
$\{d_1,...,d_s\}$ be another set of generators for $M_{\FF}$, and define in an analogous way the map $ord_{\FF,Y} '$
corresponding to this set. By assumption there exist some elements $a_{i,j}$ of $\OX$ such that:
$$
\der_i= \sum_{j=1} ^s a_{i,j} d_j
$$
By Leibniz rule, it is easy to check via an induction on $|I|$ that there exist some elements $a_{I,J}$ in $\OX$ such
that:
$$
\der_I = \sum_{|J|\leq |I|} a_{I,J} d_J
$$
Let $f$ be an element of $\OX$ such that $ord_{\FF,Y}(f)=n$. Then there exists an index $I$ of length $n$ such that:
$$
\der_I(f) = \sum_{|J|\leq n} a_{I,J} d_J(f) \; \not\in \; I_Y
$$
Since $I_Y$ is an ideal, this means there exists an index $J$ of length $\leq n$ such that $d_J(f)$ does not belong
to $I_Y$. By definition we get that $ord_{\FF,Y} '(f)\leq n=ord_{\FF,Y}(f)$ for any $f$. By symmetry we find that
$ord_{\FF,Y} '(f)=ord_{\FF,Y}(f)$ for any $f$, and the maps coincide.

\begin{prop} \label{contact}
If $Y$ is irreducible, the contact order enjoys the following properties:
\begin{itemize}
\item{$ord_{\FF,Y}(f+g)\geq \inf\{ord_{\FF,Y}(f),ord_{\FF,Y}(g)\}$ with equality if $ord_{\FF,Y}(f)\not=ord_{\FF,Y}(g)$,}
\item{$ord_{\FF,Y}(fg)=ord_{\FF,Y}(f)+ord_{\FF,Y}(g)$ for all $f,g$ in $\OX$.}
\end{itemize}
\end{prop}
{\it Proof of the first assertion}: If $ord_{\FF,Y}(f)=ord_{\FF,Y}(g)=+\infty$, then $f,g$ both belong to $I(\FF,Y)$,
$f+g$ belongs to $I(\FF,Y)$ and the result follows. So assume that $ord_{\FF,Y}(f)$ is finite and for simplicity that
$n=ord_{\FF,Y}(f)\leq ord_{\FF,Y}(g)$. For any index $I$ of length $<n$, $\der_{I}(f)$ and $\der_{I}(g)$ both belong to
$I_Y$. So $\der_{I}(f+g)$ belong to $I_Y$ for any $I$ with $|I|<n$, and $ord_{\FF,Y}(f+g)\geq n$. Therefore we have
for all $f,g$:
$$
ord_{\FF,Y}(f+g)\geq \inf\{ord_{\FF,Y}(f),ord_{\FF,Y}(g)\}
$$
Assume now that $ord_{\FF,Y}(f)<ord_{\FF,Y}(g)$. Then there exists an index $I$ of length $n$ such that $\der_{I}(f)$
does not belong to $I_Y$. Since $|I|< ord_{\FF,Y}(g)$, $\der_{I}(g)$ belongs to $I_Y$. Therefore $\der_{I}(f+g)$
does not belong to $I_Y$ and $ord_{\FF,Y}(f+g)\leq n$, so that $ord_{\FF,Y}(f+g)= n$.
\qed
For the second assertion, we will need the following lemmas. The first one is easy to get via Leibniz rule, by
an induction on the length of $I$.

\begin{lem} \label{calcul2}
Let $\der_1,...,\der_r$ a system of generators of $M_{\FF}$.
Then there exist some elements $\alpha_{I_1,I_2}$ of $\CC$, depending on $I$ and such that
for all $f,g$:
$$
\der_I (fg)= \sum_{|I_1|+|I_2|=|I|} \alpha_{I_1,I_2} \der_{I_1} (f)\der_{I_2}(g)
$$
\end{lem}

\begin{lem} \label{calcul3}
Let $f$ be an element of $\OX$ such that $ord_{\FF,Y}(f)\geq n$. Let $I=(i_1,...,i_n)$ be any index. For any
rearrangement $J=(j_1,...,j_n)$ of the $i_k$, $\der_J(f)-\der_I(f)$ belongs to $I_Y$.
\end{lem}
\dem Every rearrangement of the $i_k$ can be obtained after a composition of transpositions on two consecutive
terms. So we only need to check the lemma in the case $J=(i_1,...,i_{l+1},i_l,...,i_n)$. If we denote by
$I_1,I_2$ the indices $I_1=(i_1,...,i_{l-1})$ and $I_2=(i_{l+2},...,i_n)$, then we find:
$$
\der_J - \der_I = \der_{I_1} \circ [\der_{i_l},\der_{i_{l+1}}] \circ \der_{I_2}
$$
Since $M_{\FF}$ is stable by Lie bracket, $d=[\der_{i_l},\der_{i_{l+1}}]$ belongs to $M_{\FF}$. Then $\der_J - \der_I$
is a composite of $(n-1)$ derivations that span $M_{\FF}$. Since $ord_{\FF,Y}$ is independent of the set of generators
and $ord_{\FF,Y}(f)=n$, $\der_J (f) -\der_I(f)$ belongs to $I_Y$.
\qed
{\it Proof of the second assertion of Proposition \ref{contact}}: Let $f,g$ be a couple of elements of $\OX$. If either $f$ or $g$
has infinite contact order, then one of them belongs to $I(\FF,Y)$ and $fg$ belongs to $I(\FF,Y)$, so that
$ord_{\FF,Y}(fg)=+\infty = ord_{\FF,Y}(f) + ord_{\FF,Y}(g)$.
Assume now that $ord_{\FF,Y}(f)=n$ and $ord_{\FF,Y}(g)=m$ are finite. By lemma \ref{calcul2}, we have:
$$
\der_I (fg)= \sum_{|I_1|+|I_2|=|I|} \alpha_{I_1,I_2} \der_{I_1} (f)\der_{I_2}(g)
$$
Since $|I_1|+|I_2|<n+m$, either $|I_1|<n$ or $|I_2|<m$, and $\der_{I_1} (f)\der_{I_2}(g)$ belongs to $I_Y$. So
$\der_I (fg)$ belongs to $I_Y$ and we obtain:
$$ord_{\FF,Y}(fg)\geq n+m$$
Conversely, consider the following polynomials $P,Q$ in the indeterminates $x,t_1,...,t_r$:
$$
P(x,t_1,...,t_r)=(t_1\der_1 +...+ t_r\der_r)^n(f)(x) \quad , \quad Q(x,t_1,...,t_r)=(t_1\der_1 +...+ t_r\der_r)^m(g)(x)
$$
By lemma \ref{calcul3}, we get that $\der_I(f)\equiv \der_J(f)\; [I_Y]$ for any rearrangement $J$ of $I$ if $I$
has length $n$. Idem for $\der_I(g)$ and $\der_J(g)$ if $I$ has length $m$. Therefore in the expressions of $P,Q$,
everything happens modulo $I_Y$ as if the derivations $\der_i$ commuted. We then obtain the following expansions
modulo $I_Y$:
$$
P\equiv \sum_{i_1+...+i_r=n} \frac{n!}{i_1 ! ...i_r!} t_1 ^{i_1}...t_r ^{i_r}\der_1 ^{i_1} \circ ...\circ \der_r ^{i_r}(f) \; [I_Y]
$$
$$
Q\equiv \sum_{i_1+...+i_r=m} \frac{m!}{i_1 ! ...i_r!} t_1 ^{i_1}...t_r ^{i_r}\der_1 ^{i_1} \circ ...\circ \der_r ^{i_r}(g)\; [I_Y]
$$
Since $ord_{\FF,Y}(f)=n$ and $ord_{\FF,Y}(g)=m$, both $P$ and $Q$ have at least one coefficient that does not belong to $I_Y$
by lemma \ref{calcul3}. So neither of them belong to the ideal $I_Y[t_1,...,t_r]$, which is prime because $I_Y$ is prime. So
$PQ$ does not belong to $I_Y[t_1,...,t_r]$. If $\der=t_1 \der_1+...+t_r\der_r$, then we have by Leibniz rule:
$$
\der^{n+m}(fg)= \sum_{k=0} ^n C_{n+m} ^k \der^k(f)\der^{n+m-k} (g)
$$
Since $ord_{\FF,Y}(f)=n$ and $ord_{\FF,Y}(g)=m$, $\der^k(f)\der^{n+m-k} (g)$ belongs to $I_Y[t_1,...,t_r]$ except
for $k=n$. So $\der^{n+m}(fg)= C_{n+m} ^{n}PQ$ does not belong to $I_Y[t_1,...,t_r]$. Choose a point $(y,z_1,...,z_r)$
in $Y\times \CC^r$ such that $PQ(y,z_1,...,z_r)\not=0$ and set $d=z_1\der_1+...+z_r\der_r$. By construction we have:
$$
d^{n+m}(fg)(y)=C_{n+m} ^{n}PQ(y,z_1,...,z_r)\not=0
$$
So $d^{n+m}(fg)$ does not belong to $I_Y$ and $fg$ has contact order $\leq n+m$ with respect to the system of generators
$\{\der_1,...,\der_r,d\}$. Since the contact order does not depend on the system of generators, we find:
$$
ord_{\FF,Y}(fg)=n+m=ord_{\FF,Y}(f) +ord_{\FF,Y}(g)
$$
\qed

\begin{cor} \label{irreductible}
Let $Y$ be an irreducible subvariety of $X$. Then the ideal $I(\FF,Y)$ is prime. In particular, the minimal
invariant variety $V(\FF,Y)$ is irreducible.
\end{cor}
\dem Let $f,g$ be two elements of $\OX$ such that $fg$ belongs to $I(\FF,Y)$. Then $fg$ has infinite contact order.
By proposition \ref{contact}, either $f$ or $g$ has infinite contact order. So one of them belongs to $I(\FF,Y)$,
and this ideal is prime.
\qed

\section{Behaviour of the function $n_{\FF}$}

In this section we are going to establish theorem \ref{mesure} about the measurability of the function
$n_{\FF}$ for the Zariski topology. Recall that a function $f: X \rightarrow \mathbb{N}$ is lower
semi-continuous for the Zariski topology if the set $f^{-1}([0,r])$
is closed for any $r$. Note that such a function is continuous for the constructible topology.
We begin with the following lemma.

\begin{lem} \label{etape}
Let $F$ be a finite dimensional vector subspace of $\OX$. Then the map $\varphi_F: X \rightarrow
\mathbb{N}, \; x\mapsto dim_{\CC} \;F -dim_{\CC} \; I(\FF,x)\cap F$ is lower semi-continuous for the Zariski topology.
\end{lem}
\dem For any fixed finite-dimensional vector space $F$, consider the affine algebraic set:
$$
\Sigma_F=\left\{ (x,f) \in X\times F, \; \forall d_1,...,d_m \in M_{\FF}, \; d_1\circ...\circ d_m(f)(x)=0\right\}
$$
together with the projection $\Pi: \Sigma_F \longrightarrow X, \; (x,f)\longmapsto x$. Since $\Sigma_F$ is affine,
there exists a finite collection of linear operators $\Delta_1,...,\Delta_r$, obtained by composition of elements of
$M_{\FF}$,
such that:
$$
\Sigma_F=\left\{ (x,f) \in X\times F, \;  \Delta_1(f)(x)=...=\Delta_r(f)(x)=0\right\}
$$
By lemma \ref{autre}, the fibre $\Pi^{-1}(x)$ is isomorphic to $I(\FF,x)\cap F$ for any point $x$ of $X$. Since every
$\Delta_i$ is linear, $\Delta_i$ can be considered as a linear form on $F$ with coefficients in $\OX$. So
the map $\Delta=(\Delta_1,...,\Delta_r)$ is represented by a matrix with entries in $\OX$. We therefore have the
equivalence:
$$
f\in I(\FF,x)\cap F \quad \Longleftrightarrow \quad f \in ker \;  \Delta(x)
$$
By the rank theorem, we have $\varphi_F(x)=rk \; \Delta(x)$. But the rank of this matrix is a lower semi-continuous function
because it is given as the maximal size of the minors of $\Delta$ that do not vanish at $x$. Therefore $\varphi_F$
is lower semi-continuous for the Zariski topology.
\qed
{\it Proof of theorem \ref{mesure}}: Since $X$ is affine, we may assume that $X$ is embedded in $\CC^k$ for some $k$. We provide
$\CC[x_1,...,x_k]$ with the filtration $\{F_n\}$ given by the polynomials of homogeneous degree $\leq n$. By Hilbert-Samuel
theorem (see \cite{Ei}), for any ideal $I$ of $\CC[x_1,...,x_k]$, the function:
$$
h_I(n)=dim_{\CC} \; F_n - dim_{\CC} \; I\cap F_n
$$
is equal to a polynomial for $n$ large enough, and the degree $p$ of this polynomial coincides with the dimension of the variety
$V(I)$. It is therefore easy to show that:
$$
p= \lim_{n\to +\infty} \frac{\log(h_I(n))}{n}
$$
Let $\Pi: \CC[x_1,...,x_k]\rightarrow \OX$ be the morphism induced by the inclusion $X\hookrightarrow \CC^k$, and set
$\widetilde{F_n}=\Pi(F_n)$. For any ideal $I$ of $\OX$, consider the function:
$$
\widetilde{h_I}(n)=dim_{\CC} \; \widetilde{F_n} - dim_{\CC} \; I\cap \widetilde{F_n}
$$
Since $\Pi$ is onto, we have $\widetilde{h_I}(n)=h_{\Pi^{-1}(I)}(n)$, so that $\widetilde{h_I}(n)$ coincides for $n$ large
enough with a polynomial of degree $p$ equal to the dimension of $V(I)$. With the notation of lemma \ref{etape}, we obtain
for $I=I(\FF,x)$:
$$
p=n_{\FF}(x)= \lim_{n\to +\infty} \frac{\log(\widetilde{h_I}(n))}{n}=\lim_{n\to +\infty} \frac{\log(\varphi_{\widetilde{F_n}}(x))}{n}
$$
By lemma \ref{etape}, every $\varphi_{\widetilde{F_n}}$ is lower semi-continuous for the Zariski topology, hence measurable.
Since a pointwise limit of measurable functions is measurable, the function $n_{\FF}$ is measurable for the Zariski topology.
Moreover since $\varphi_{\widetilde{F_n}}$ is lower semi-continuous, there exist a real number $r_n$ and an open set $U_n$ on
$X$ such that:
\begin{itemize}
\item{$\displaystyle \frac{\log(\varphi_{\widetilde{F_n}}(x))}{n}\leq r_n$ for any $x$ in $X$,}
\item{$\displaystyle \frac{\log(\varphi_{\widetilde{F_n}}(x))}{n}=r_n$ for any $x$ in $U_n$.}
\end{itemize}
Denote by $U$ the intersection of all $U_n$. Since this intersection is not empty, there exists an $x$ in $X$ for which
$\log(\varphi_{\widetilde{F_n}}(x))/n=r_n$ for any $n$, so that $r_n$ converges to a limit $p$. By passing to the limit,
we obtain that:
\begin{itemize}
\item{$n_{\FF}(x)\leq p$ for any $x$ in $X$,}
\item{$n_{\FF}(x)= p$ for any $x$ in $U$.}
\end{itemize}
Note that $p$ has to be an integer. The theorem is proved.
\qed

\section{The family of minimal invariant varieties} \label{set}

In this section, we are going to study the set of minimal invariant varieties associated to the points of $X$. The
result we will get will be the first step towards the proof of theorem \ref{fibration}. Let $M$ be the following set:
$$
M=\left\{(x,y) \in X\times X, \; y \in V(\FF,x)\right\}
$$
together with the projection $\Pi: M\longrightarrow X, (x,y)\longmapsto x$. Note that for any $x$, the preimage $\Pi^{-1}(x)$
is isomorphic to $V(\FF,x)$, so that the couple $(M,\Pi)$ parametrizes the set of all minimal invariant varieties. Our purpose is
to show that: 

\begin{prop} \label{ferm}
The Zariski closure $\overline{M}$ is an irreducible affine set of dimension $dim\; X +p$, where $p$ is the maximum of the function
$n_{\FF}$. Moreover, for every very generic point $x$ in $X$, $\overline{M} \cap \Pi ^{-1}(x)$ is equal to $\{x\}\times V(\FF,x)$.
\end{prop}
The proof of this proposition is a direct consequence of the following lemmas.

\begin{lem}
The Zariski closure $\overline{M}$ is irreducible.
\end{lem}
\dem For any $\der_i$ in $\FF$, consider the new $\CC$-derivation $\Delta_i$ on ${\cal{O}}_{X\times X}=\OX \otimes_{\CC} \OX$
given by the following formula:
$$
\forall f,g \in \OX, \quad \Delta_i(f(x)\otimes g(y))=f(x)\otimes \der_i(g)(y)
$$
It is easy to check that $\Delta_i$ is a well-defined derivation. Denote by ${\cal{G}}$ the collection of the $\Delta_i$, by $D$
the diagonal $\{(x,x), \; x \in X\}$ in $X\times X$ and set $M_0= V({\cal{G}},D)$. By corollary \ref{irreductible}, $M_0$ is
irreducible.
We are going to prove that $\overline{M}=M_0$.

First let us check that $M_0 \subseteq \overline{M}$. Let $f$ be a regular function on $X\times X$ that vanishes on $\overline{M}$.
Then $f(x,y)=0$ for any couple $(x,y)$ where $y$ belongs to $V(\FF,x)$. If $\varphi_t(y)$ is the flow of $\der_i$ at $y$, then
$\psi_t(x,y)=(x,\varphi_t(y))$ is the flow of $\Delta_i$ at $(x,y)$. Since $y$ lies in $V(\FF,x)$, $\varphi_t(y)$ belongs to
$V(\FF,x)$ for any small value of $t$, and we obtain:
$$
f(\psi_t(x,y))=f(x,\varphi_t(y))=0
$$
By derivation with respect to $t$, we get that $\Delta_i(f)(x,y)=0$ for any $(x,y)$ in $M$. So $\Delta_i(f)$ vanishes along
$\overline{M}$, and the ideal $I(\overline{M})$ is stable by the family ${\cal{G}}$. Since it is contained in $I(D)$, we have the
inclusion:
$$
I(\overline{M})\subseteq I({\cal{G}},D)
$$
which implies that $M_0 \subseteq \overline{M}$.

Second let us show that $\overline{M}\subseteq M_0$. Let $f$ be a regular function that vanishes along $M_0$. Fix $x$ in $X$
and consider the function $f_x(y)=f(x,y)$ on $X$. Then for any $\Delta_1,...,\Delta_n$ in ${\cal{G}}$, we have:
$$
\Delta_1 \circ ...\circ \Delta_n(f)(x,y)= \der_1\circ ...\circ \der_n(f_x)(y)
$$
Since $M_0=V({\cal{G}},D)$, $D$ is contained in $M_0$ and $f_x(x)=0$. So $f(x,x)=0$ and for any $\der_1,...,\der_n$ in $\FF$ and
any $x$ in $X$, we get that:
$$
\der_1\circ ...\circ \der_n(f_x)(x)=0
$$
In particular, $f_x$ belongs to $I(\FF,x)$ and $f_x$ vanishes along $V(\FF,x)$. Thus $f$ vanishes on $\{x\}\times
V(\FF,x)=\Pi^{-1}(x)$ for any $x$ in $X$. This implies that $f$ is equal to zero on $M$ and on $\overline{M}$, so that
$I({\cal{G}},D)\subseteq I(\overline{M})$. As a consequence, we find $\overline{M}\subseteq M_0$ and the result follows.
\qed

\begin{lem}
The variety $\overline{M}$ has dimension $\geq dim \; X+p$.
\end{lem}
\dem Consider the projection $\Pi: \overline{M} \rightarrow X, \; (x,y)\mapsto x$. Since $M$ contains the diagonal $D$,
the map $\Pi$ is onto. By the theorem on the dimension of fibres, there exists a non-empty Zariski open set $U$ in $X$
such that:
$$
\forall x \in U, \quad dim \; \overline{M} = dim\; X + dim \; \Pi^{-1}(x) \cap \overline{M}
$$
By theorem \ref{mesure}, there exists a countable intersection $\theta$ of Zariski open sets in $X$ such that
$n_{\FF}(x)=p$ for all $x$ in $X$. In particular, $U \cap \theta$ is non-empty. For any $x$ in
$U \cap \theta$, $\Pi^{-1}(x) \cap \overline{M}$ contains the variety $V(\FF,x)$ whose dimension
is $p$, and this yields:
$$
dim \; \overline{M} \geq  dim\; X + p
$$
\qed

\begin{lem} \label{note}
The variety $\overline{M}$ has dimension $\leq dim \; X+p$.
\end{lem}
\dem Let $\{F_n\}$ be a filtration of $\OX$ by finite-dimensional $\CC$-vector spaces, and set:
$$
M_n =\left\{ (x,y) \in X\times X, \; \forall f \in I(\FF,x)\cap F_n, \; f(y)=0 \right\}
$$
The sequence $\{M_n\}$ is decreasing for the inclusion, and $M=\cap_{n \in \mathbb{N}}\;  M_n$. Moreover
every $M_n$ is constructible for the Zariski topology by Chevalley's theorem (see \cite{Ei}).
Indeed its complement in $X\times X$ is the
image of the constructible set:
$$
\Sigma_n =\left\{ (x,y,f) \in X\times X \times F_n , \; \forall \der_1,...,\der_k \in \FF, \; \der_1 \circ ...\circ
\der_k(f)(y)=0 \; {\rm{and}} \; f(y)\not=0\right\}
$$
under the projection $(x,y,f)\mapsto (x,y)$. Since $D$ is contained in every $M_n$, the projection $\Pi: M_n \rightarrow
X$ is onto. By the theorem on the dimension of fibres applied to the irreducible components of $\overline{M_n}$, there
exists a non-empty Zariski open set $U_n$ in $X$ such that:
$$
\forall x \in U_n, \quad dim \; M_n \leq dim\; X + dim \; \Pi^{-1}(x) \cap M_n
$$
Since $\overline{M}\subseteq \overline{M_n}$ for any $n$, and $\Pi^{-1}(x) \cap M_n\simeq V(I(\FF,x)\cap F_n)$, we obtain:
$$
\forall x \in U_n, \quad dim \; \overline{M} \leq dim\; X + dim \; V(I(\FF,x)\cap F_n)
$$
Since every $U_n$ is open, the intersection $\theta'=\cap_{n\in \mathbb{N}} \; U_n$ is non-empty. Let $\theta$ be an intersection
of Zariski open sets of $X$ such that $n_{\FF}(x)=p$ for any $x$ of $\theta$. For any fixed $x$ in
$\theta \cap \theta'$, we have:
$$
\forall n\in \mathbb{N}, \quad dim \; \overline{M} \leq dim\; X + dim \; V(I(\FF,x)\cap F_n)
$$
Since $\OX$ is noetherian, there exists an order $n_0$ such that $I(\FF,x)$ is generated by $I(\FF,x)\cap F_n$ for any
$n\geq n_0$. In this context, $V(\FF,x)=V(I(\FF,x)\cap F_n)$ for all $n\geq n_0$, and $V(\FF,x)$ has dimension $p$,
which implies that:
$$
dim \; \overline{M} \leq dim\; X + p
$$
\qed

\begin{lem}
For every very generic point $x$ in $X$, $\overline{M} \cap \Pi ^{-1}(x)$ is equal to $\{x\}\times V(\FF,x)$.
\end{lem}
\dem Consider the constructible sets $M_n$ introduced in lemma \ref{note}. By construction their intersection is equal to $M$.
The $\{\overline{M_n}\}$ form a decreasing sequence which converges to $\overline{M}$. Since these are algebraic sets, there exists
an index $n_0$ such that for any $n\geq n_0$, we have $\overline{M_n}=\overline{M}$. We consider the sequence $\{M_n\}_{n\geq n_0}$
and denote by $G_n$ the Zariski closure of $\overline{M} - M_n$. By the theorem on the dimension of fibres, there exists a Zariski
open set $V_n$ on $X$ such that for any $x$ in $V_n$, either
$\Pi^{-1}(x)\cap G_n$ is empty or has dimension $<p$. Since $\Pi^{-1}(x)\cap M=\{x\}\times V(\FF,x)$ for any $x$ in $X$, we have the following
decomposition:
$$
\Pi^{-1}(x)\cap \overline{M}= \{x\}\times V(\FF,x) \cup \cup_{n\geq n_0} \Pi^{-1}(x)\cap G_n
$$
For all $x$ in $\theta=\cap V_n$, the set $\Pi^{-1}(x)\cap G_n$ has dimension $<p$ for any $n\geq n_0$, hence its Hausdorff dimension
is no greater than $(2p-2)$ (see \cite{Ch}). Consequently the countable union $\cup_{n\geq n_0} \Pi^{-1}(x)\cap G_n$ has an Hausdorff dimension
$<2p$. Let $H_{i,x}$ be the irreducible components of $\Pi^{-1}(x)\cap \overline{M}$ distinct from $\{x\}\times V(\FF,x)$. These $H_{i,x}$
are covered
by the union $\cup_{n\geq n_0} \Pi^{-1}(x)\cap G_n$, hence their Hausdorff dimension does not exceed $(2p-2)$. Therefore the Krull
dimension of $H_{i,x}$ is strictly less than $p$ for any $i$ and any $x$ in $\theta$. If $H_x$ denotes the union of the $H_{i,x}$, then
we have for any $x$ in $\theta$:
$$
\Pi^{-1}(x)\cap \overline{M}= \{x\}\times V(\FF,x) \cup H_x \quad \mbox{and} \quad dim\; H_x<p
$$
Now by Stein factorization theorem (see \cite{Ha}), the map $\Pi: \overline{M}\rightarrow X$ is a composite of a quasi-finite map
with a map whose generic fibres are irreducible. In particular $\Pi^{-1}(x)\cap \overline{M}$ is equidimensionnal of dimension $p$
for generic $x$ in $X$. Therefore the variety $H_x$ should be contained in $\{x\}\times V(\FF,x)$, and we have for any $x$ in $\theta$:
$$
\Pi^{-1}(x)\cap \overline{M}= \{x\}\times V(\FF,x) 
$$
\qed

\section{Proof of theorem \ref{fibration}}

Let $X$ be an irreducible affine variety over $\CC$ of dimension $n$, endowed with an algebraic foliation $\FF$. Let $p$ be the integer
given by theorem \ref{mesure}. In this section we will establish theorem \ref{fibration}. We begin with
a few lemmas.

\begin{lem}
Let $F:X\rightarrow Y$ be a dominant morphism of irreducible affine varieties. Then for any Zariski open set $U$ in $X$,
$F(U)$ is dense in $Y$.
\end{lem}
\dem Suppose on the contrary that $F(U)$ is not dense in $Y$. Then there exists a non-zero regular function $f$ on $Y$
that vanishes along $\overline{F(U)}$. The function $f\circ F$ vanishes on $U$, hence on $X$ by density. So $F(X)$ is
contained in $f^{-1}(0)$, which is impossible since this set is dense in $Y$.
\qed

\begin{lem} \label{prep}
Let $\overline{M}$ be the variety defined in section \ref{set}. Then there exists an irreducible variety $H$
in $X$ such that $\overline{M} \cap X\times H$ has dimension $n$ and the morphism $\Pi: \overline{M} \cap X\times H
\rightarrow X$ induced by the projection is dominant. 
\end{lem}
\dem Let $(x,y)$ be a smooth point of $\overline{M}$ such that $x$ is a smooth point of $X$. By the generic smoothness
theorem, we may assume that $d\Pi_{(x,y)}$ is onto. Consider the second projection $\Psi(x,y)=y$. Since the map $(\Pi,\Psi)$
defines an embedding of $\overline{M}$ into $X\times X$, and $d\Pi_{(x,y)}$ is onto, there exist some regular functions
$g_1,...,g_p$ on $X$ such that $(d\Pi_{(x,y)},{dg_1}_{(y)},...,{dg_p}_{(y)})$ is an isomorphism from $T_{(x,y)} \overline{M}$ to
$T_x X \oplus \CC ^p$.

Let $G: \overline{M}\rightarrow \CC^p$ be the map $(g_1,...,g_p)$, and denote by $E$ the set of points $(x,y)$ in
$\overline{M}$
where either $\overline{M}$ is singular
or $(\Pi,G)$ is not submersive. By construction $E$ is a closed set distinct from $\overline{M}$.
Since $dG_{(y)}$ has rank $p$ on $T_{(x,y)} \overline{M}$, the map $G:\overline{M}\rightarrow \CC^p$
is dominant. So its generic fibres have dimension $n$. Fix a fibre $G^{-1}(z)$ of dimension $n$
that is not contained in $E$. Then there exists a smooth point $(x,y)$ in $G^{-1}(z)$ such that
$d(\Pi,G)_{(x,y)}$ is onto. The morphism $\Pi:G^{-1}(z)\rightarrow X$ is a submersion at $(x,y)$,
hence it is dominant. Moreover $G^{-1}(z)$ is of the form $X\times F^{-1}(z)\cap \overline{M}$, where $F:X
\rightarrow \CC^p$ is the map $(g_1,...,g_p)$.

Choose an irreducible component $H$ of $F^{-1}(z)$ such that $\Pi: X\times H \cap \overline{M}\rightarrow X$ is dominant.
By construction $X\times H \cap \overline{M}$ has dimension $\leq n$. Since the latter map is dominant,
its dimension is exactly equal to $n$. 
\qed
{\it Proof of theorem \ref{fibration}}: Let $H$ be an irreducible variety of codimension $p$ in $X$ satisfying the
conditions of lemma \ref{prep}. Denote by
$N$ the union of irreducible components of $\overline{M} \cap X\times H$ that are mapped dominantly on $X$ by $\Pi$.
By construction $N$ has dimension $dim\; X$ and the morphism $\Pi: N\rightarrow X$ is quasi-finite.
So there exists an open set $U$ in $X$ such that:
$$
\widetilde{\Pi}: \Pi^{-1}(U)\cap N \longrightarrow U
$$
is a finite unramified morphism. Let $r$ be the degree of this map. For any point $x$ in $U$, there exist $r$ points
$y_1,...,y_r$ in $H$ such that ${\widetilde{\Pi}}^{-1}(x)=\{y_1,...,y_r\}$. Let $\mathfrak{S}_r$ act on $H^r$ by
permutation
of the coordinates, i.e $\sigma.(y_1,...,y_r)=(y_{\sigma(1)},...,y_{\sigma(r)})$. Since this action is algebraic
and $\mathfrak{S}_r$ is finite, the algebraic quotient $H^r //\mathfrak{S}_r$ exists and is an irreducible affine
variety (see \cite{Mu}). Let $Q: H\rightarrow H^r //\mathfrak{S}_r$ be the corresponding quotient morphism.
Consider the mapping:
$$
\varphi: U \longrightarrow H^r //\mathfrak{S}_r , \quad x \longmapsto Q(y_1,...,y_r)
$$
Note that its graph is constructible in $U\times H^r //\mathfrak{S}_r$. Indeed it is given by the set:
$$
\Sigma=\left\{ (x,y'), \; \exists (y_1,...,y_r) \in H^r, \;\forall i\not=j, \; y_i\not=y_j, \; (x,y_i) \in \overline{M} \; 
{\rm{and}} \; Q(y_1,...,y_r)=y' \right\}
$$
By Serre's theorem (see \cite{Lo}), $\varphi$ is a rational map on $U$. Since $\widetilde{\Pi}$ is unramified, $\varphi$ is also
holomorphic on $U$, hence it is regular on $U$. Denote by $Y$ the Zariski closure of $\varphi(U)$ in
$H^r //\mathfrak{S}_r$. Since $U$ is irreducible, $Y$ is itself irreducible.

By construction, for any $x$ in $U$, $\{x\}\times \varphi^{-1}(\varphi(x))$ is equal to $\Pi^{-1}(x)\cap \overline{M}$. 
For every very generic point $x$ in $X$, $\Pi^{-1}(x)\cap \overline{M}$ corresponds to $\{x\}\times V(\FF,x)$
by proposition \ref{ferm}.
So $\varphi^{-1}(\varphi(x))=V(\FF,x)$ for every generic point $x$ in $X$, hence it has dimension $p$. By the
theorem on the dimension of fibres, $Y$ has dimension $(n - p)$.

Since $\varphi^{-1}(\varphi(x))=V(\FF,x)$ for every generic point $x$ in $X$, this fibre is tangent to the
foliation $\FF$. Since tangency is a closed condition, all the fibres of $\varphi$ are tangent to $\FF$. Let
$f$ be a rational function on $Y$. In the neighborhood of any smooth point $x$ where $\FF$ is regular and $f\circ \varphi$ is
well-defined, the function $f\circ \varphi$ is constant on the leaves of $\FF$. So $f\circ \varphi$ is a rational first
integral of $\FF$. Via the morphism $\varphi^*$ induced by $\varphi$, $K_{\FF}$ is clearly isomorphic to $\CC(Y)$
which has transcendence degree $(n-p)$ over $\CC$.
\qed

\section{An example}

In this last section, we introduce an example that illustrates both theorems \ref{mesure} and \ref{fibration}. Consider the
affine space $\CC^4$ with coordinates $(u,v,x,y)$, and the algebraic foliation $\FF$ induced by the vector field:
$$
\der=ux\frac{\der}{\der x} + vy \frac{\der}{\der y}
$$
For any $(\lambda,\mu)$ in $\CC^2$, the plane $V(u-\lambda,v-\mu)$ is tangent to $\FF$. Denote by $\der_{\lambda,\mu}$ the
restriction of $\der$ to that plane parametrized by $(x,y)$. Then two cases may occur:
\begin{itemize}
\item{If $[\lambda;\mu]$ does not belong to $\mathbb{P}^1(\mathbb{Q})$, then $\der_{\lambda,\mu}$ has no rational first integrals.
The only algebraic curves tangent to $\der_{\lambda,\mu}$ are the lines $x=0$ and $y=0$. There is only one singular point, namely
$(0,0)$.}
\item{If $[\lambda;\mu]$ belongs to $\mathbb{P}^1(\mathbb{Q})$, choose a couple of coprime integers $(p,q)\not=(0,0)$ such that
$p\lambda+q\mu=0$. The function $f(x,y)=x^py^q$ is a rational first integral for $\der_{\lambda,\mu}$. The algebraic curves tangent
to $\der_{\lambda,\mu}$ are the lines $x=0$, $y=0$ and the fibres $f^{-1}(z)$ for $z\not=0$. There is only one singular point, namely
$(0,0)$.}
\end{itemize}
From those two cases, we can get the following values for the function $n_{\FF}$:

\begin{itemize}
\item{$n_{\FF}(u,v,x,y)=2$ if $[\lambda;\mu] \not\in\mathbb{P}^1(\mathbb{Q})$ and $xy\not=0$,}
\item{$n_{\FF}(u,v,x,y)=0$ if $x=y=0$,}
\item{$n_{\FF}(u,v,x,y)=1$ otherwise.}
\end{itemize}
In particular, this function is measurable but not constructible for the Zariski topology, as can be easily seen
from its fibre $n_{\FF} ^{-1}(2)$. Moreover since $p=2$, its
field $K_{\FF}$ has transcendence degree 2 over $\CC$. In fact it is easy to check that $K_{\FF}=\CC(u,v)$.


\begin{thebibliography}{G-M}

\bibitem[Bru]{Bru} M.Brunella, M.Nicolau {\it Sur les hypersurfaces solutions des \'equations de Pfaff}, note
aux C.R.A.S t.329 s\'erie I, pp. 793-795 (1999).

\bibitem[Ch]{Ch} E.M.Chirka {\it Complex analytic sets}, Kluwer Academic Publishers Group, Dordrecht (1989).

\bibitem[Ei]{Ei} D.Eisenbud {\it Commutative Algebra with a view toward Algebraic
Geometry}, Springer Verlag New York (1995).

\bibitem[Gh]{Gh} E.Ghys {\it A propos d'un th\'eor\`eme de Jouanolou concernant les feuilles ferm\'ees des feuilletages
holomorphes}, Rend. Circ. Mat. di Parlermo (1999).

\bibitem[G-M]{G-M} X.Gomez-Mont {\it Integrals for holomorphic foliations with singularities having all leaves compact},
Ann. Inst. Fourier (Grenoble) 39, $n^o$ 2, pp. 451-458 (1989).

\bibitem[Ha]{Ha} R.Hatshorne {\it Algebraic geometry}, Springer Verlag New York (1977).

\bibitem[Jou]{Jou} J-P.Jouanolou {\it Hypersurfaces solutions d'une \'equation de Pfaff analytique}, Math.Ann.
239, pp. 239-245 (1978).

\bibitem[Ka]{Ka} I.Kaplansky {\it An introduction to Differential Algebra}, publications de l'institut
de math\'ematiques de Nancago, Hermann (1957).

\bibitem[Lo]{Lo} S.Lojasiewicz {\it Introduction to complex analytic geometry}, Birkh\"auser Verlag, Basel (1991).

\bibitem[Mu]{Mu} D.Mumford, J.Fogarty, F.Kirwan {\it Geometric invariant theory}, third edition Springer Verlag, Berlin (1994).

\bibitem[Pe]{Pe} J.V.Pereira {\it Vector fields, invariant varieties and linear systems}, Ann.Inst.Fourier (Grenoble) 51, $n^0$ 5,
pp.1385-1405 (2001).


\end{thebibliography}
\end{document}